\documentclass[12pt]{amsart}
\usepackage{amssymb,latexsym,amsthm}
\usepackage{graphicx}

\newcommand{\mfA}{\mathfrak{A}}
\newcommand{\mfB}{\mathfrak{B}}
\newcommand{\rad}{\mathrm{rad}}

\newtheorem{theorem}{Theorem}[section]
\newtheorem{lemma}[theorem]{Lemma}
\newtheorem{cor}[theorem]{Corollary}

\theoremstyle{definition}

\newtheorem{example}[theorem]{Example}

\theoremstyle{remark}

\numberwithin{equation}{section}

\begin{document}

\baselineskip=17pt

\title[]
{Linear extensions of isometries between groups of  
invertible elements in Banach algebras}

\author{Osamu~Hatori}
\address{Department of Mathematics, Faculty of Science, 
Niigata University, Niigata 950-2181 Japan}
\curraddr{}
\email{hatori@math.sc.niigata-u.ac.jp}

\thanks{The author was partly 
supported by the Grants-in-Aid for Scientific 
Research, The 
Ministry of Education, Science, Sports and Culture, Japan.}



\maketitle

\begin{abstract}
We show that if $T$ is an isometry (as metric spaces) 
from an 
open subgroup of 
the invertible group $A^{-1}$ of a unital Banach algebra $A$
onto an 
open subgroup of 
the invertible group $B^{-1}$ of a unital Banach algebra $B$, 
then $T$ is 
extended to a real-linear isometry up to translation between these 
Banach algebras. 
We consider multiplicativity or unti-multiplicativity 
of the isometry. Note that a unital linear isometry between unital semisimple 
commutative Banach algebra need be multiplicative. 
On the other hand, we show that 
 if $A$ is commutative and $A$ or $B$ are semisimple, 
then 
$(T(e_A))^{-1}T$ is extended to a isometrical real algebra isomorphism
from $A$ onto $B$. In particular, 
$A^{-1}$ is isometric as a metric space to $B^{-1}$ if and only if 
they are isometrically isomorphic to each other as metrizable groups if and 
only if 
$A$ is 
isometrically isomorphic to $B$ as a real Banach algebra; 
it is compared by the example of 
\.Zelazko concerning on non-isomorphic Banach algebras with 
homeomorphically isomorphic invertible groups. 
Maps between standard operator algebras are also investigated.
\end{abstract}
\section{Introduction}
A long tradition of inquiry seeks sufficient sets of 
conditions on Banach algebras and related subsets, in order 
that isometries on the algebras or subsets 
are automatically (or extended to be) linear, multiplicative, or 
Jordan.
In this paper we investigate 
isometries between groups of invertible elements in 
unital complex Banach algebras.

The Mazur-Ulam theorem asserts that every isometry from a normed 
space onto a normed space is real linear up to translation.
Suppose that 
a map $T$ is defined from an open set of a Banach space 
onto an open set of a Banach space. Even if $T$ is isometric, 
$T$ need not be extended a real linear map up to 
translation between the 
Banach spaces (see Example \ref{hoo}).
We show in the section \ref{lin} 
that an isometry $T$ from an open subgroup of the invertible 
group of a unital Banach algebra $A$ onto an open subgroup of 
the invertible group of a unital Banach algebra $B$ is 
extended to a real-linear isometry up to 
translation of $A$ onto $B$. 

There exist unital complex Banach algebras $A$ and $B$ which are 
not isomorphic as real algebras while the invertible groups 
$A^{-1}$ and $B^{-1}$ are homeomorphically isomorphic to 
each other as topological groups (cf. \cite[Remark 1.7.8.]{z2}).
In the section \ref{met} 
we show that a unital complex Banach algebra 
$A$ is isometrically isomorphic to a unital complex 
Banach algebra $B$ as a real Banach algebra if and only if 
$A^{-1}$ is isometrically isomorphic to $B^{-1}$ as a 
metrizable group if and only if the principal component of 
$A^{-1}$ is isometrically isomorphic to the 
principal component of $B^{-1}$ as a metrizable group. 

In the section \ref{mul} we consider multiplicativity of 
the extended linear map obtained in Theorem \ref{main}.
Suppose that the given algebra $A$ in Theorem \ref{main} is a 
uniform algebra. Then the extended map $\tilde T_0$ is
multiplicative if $T$ is unital, $T(ia)=iT(a)$ for 
every $a\in \mfA$, and $u_0=0$ by a 
theorem of Jarosz \cite[Corollary 2]{ja1}, 
which is a generalization of the Banach-Stone theorem. 
On the other hand a unital isometry of unital semisimple 
commutative 
Banach algebras need not be multiplicative even if 
it is surjective and complex-linear; a Banach-Stone 
theorem does not hold for general 
unital semisimple commutative Banach algebras. 
To prove multiplicativity of $\tilde T_0$ we need 
alternative way.
In \cite{hatori} we proved by applying a theory of 
spectrum preserving map that the invertible 
group ${\mathcal A}^{-1}$ of a unital {\it semisimple commutative} 
Banach algebra is isometric to the 
invertible group $B^{-1}$ of a unital Banach algebra 
as a metric space if and 
only if ${\mathcal A}^{-1}$ is isometrically isomorphic to 
$B^{-1}$ as a metrizable group if and only if 
${\mathcal A}$ is isometrically isomorphic to $B$ as 
a {\it real} Banach algebra. 
In the section \ref{mul} we show a comparison result 
(Corollary \ref{comsem}) to the above. 
We also consider the case for closed standard 
operator algebras.

Throughout the paper we denote the invertible group, 
the group of all the invertible elements in a given Banach 
algebra $A$, by $A^{-1}$. The principal component of 
$A^{-1}$ is denoted by $A^{-1}_e$. 
Note that every open subgroup of the invertible 
group contains the principal component for every unital 
Banach algebras. 
For $a\in A$, $\sigma (a)$ denotes 
the spectrum of $a$ and the spectral radius of $a$ 
is $r (a)$. The unit element in $A$ is denoted by $e_A$.
The (Jacobson) radical of $A$ is $\rad (A)$.
\section{Lemmata}
We begin by showing 
a preliminary lemma. We prove it by making use of an idea of 
V\"ais\"al\"a \cite{v}
\begin{lemma}\label{vvl}
Let $c\in B$ and a map $\psi:B\to B$ be defined as 
$\psi (z)=2c-z$ Suppose that $L$ is a non-empty 
bounded subset of $B$ such that $c\in L$ and $\psi (L)=L$. 
If ${\mathcal T}$is a 
surjective isometry from $L$ onto itself.
Then ${\mathcal T}(c)=c$.
\end{lemma}
\begin{proof}
Let $W$ be the set of all surjective isometries from $L$ onto 
itself. Note that $W$ is not empty since the identity function 
is in $W$. Let
\[
\lambda =\sup\{\|g(c)-c\|:g\in W\}.
\]
Since $L$ is non-empty and bounded $\lambda$ is not infty. 
We will show that $\lambda=0$. It will force that 
$g(c)=c$ for every surjective isometry from $L$ onto 
itself. 
Suppose that 
$g\in W$. Let $g^*=g^{-1}\circ \psi \circ g$. Then 
$g^*\in W$. Hence 
\begin{multline*}
\lambda\ge \|g^*(c)-c\|=\|g^{-1}\circ \psi\circ g(c)
-c\|\\
=\|\psi \circ g(c)-g(c)\|=2\|g(c)-c\|.
\end{multline*}
It follows that $\lambda \ge 2\lambda$ since $g$ can be chosen 
arbitrary, hence $\lambda =0$.
\end{proof}
A real linear space with a metric 
$d(\cdot,\cdot)$ satisfying $d(a+u,b+u)=d(a,b)$ 
for all $a,b,u$, and for which addition and 
scalar multiplication are jointly continuous is 
called {\it a metric real linear space}.
\begin{lemma}\label{lmu}
Let ${\mathcal B}_1$ be a real normed space and 
${\mathcal B}_2$ a metric real linear space with a 
metric $d(\cdot,\cdot)$. Suppose that 
$U_1$ and $U_2$ are non-empty open subsets of 
${\mathcal B}_1$ and 
${\mathcal B}_2$ respectively. Suppose that ${\mathcal T}$ is 
a surjective isometry 
\rm{(}$d({\mathcal T}(a),{\mathcal T}(b))=
\|a-b\|$ for 
every $a,b\in U_1$\rm{)} from $U_1$ onto $U_2$ 
and $f,g\in U_1$. If $f$ and $g$  satisfy 
the equation $(1-r)f+rg\in U_1$ for every 
$r$ with $0\le r \le 1$, then the equality
\[
{\mathcal T}(\frac{f+g}{2})=\frac{{\mathcal T}(f)+{\mathcal T}(g)}{2}
\]
holds.
\end{lemma}
\begin{proof}
Let $h,h'\in U_1$. Suppose that there exists $\varepsilon >0$ 
which satisfies that 
$\frac{\|h-h'\|}{2}<\varepsilon$, and 
\[
\{a\in {\mathcal B}_1:\|a-h\|<\varepsilon,\,\,\|a-h'\|<\varepsilon \}\subset U_1,
\]
\[
\{u\in {\mathcal B}_2:d(u,{\mathcal T}(h))<\varepsilon,\,\,
d(u,{\mathcal T}(h'))<\varepsilon \}\subset U_2.
\]
We will show that 
${\mathcal T}(\frac{h+h'}{2})=\frac{{\mathcal T}(h)+{\mathcal T}(h')}{2}$.
Set 
$r=\frac{\|h-h'\|}{2}$ and let
\[
L_1=
\{a\in {\mathcal B}_1:\|a-h\|=r=\|a-h'\|\},
\]
\[
L_2=
\{u\in {\mathcal B}_2:d(u,{\mathcal T}(h))=r=d(u,{\mathcal T}(h'))\}.
\]
Set also $c_1=\frac{h+h'}{2}$ and 
$c_2=\frac{{\mathcal T}(h)+
{\mathcal T}(h')}{2}$. Then we have ${\mathcal T}(L_1)=L_2$, 
$c_1\in L_1 \subset U_1$, and $L_2 \subset U_2$. 
Let 
\[
\psi_2(y)=2c_2-y \quad (y\in {\mathcal B}_2).
\]
Then $\psi_2$ is an isometry and $\psi_2(L_2)=L_2$. 
Let $Q={\mathcal T}^{-1}\circ\psi_2\circ{\mathcal T}$. 
Then $Q$ is well-defined and is a surjective isometry 
from  $L_1$ onto itself. 
Then by Lemma \ref{vvl} $Q(c_1)=c_1$. Henceforth 
${\mathcal T}(c_1)=c_2$.

We assume that $f$ and $g$ are as described. Let
\[
K=\{(1-r)f+rg:0\le r\le 1\}.
\]
Since $K$ and ${\mathcal T}(K)$ are compact, there is 
$\varepsilon >0$ with
\[
\inf\{\|a-b\|:a\in K,\,b\in {\mathcal B}_1\setminus U_1\}
>\varepsilon, 
\]
\[
\inf\{d(u,v):u\in {\mathcal T}(K),\,
v\in {\mathcal B}_2\setminus U_2\}>
\varepsilon.
\]
Then for every $h\in K$ we have 
\[
\{a\in {\mathcal B}_1:\|a-h\|<\varepsilon\}\subset U_1
\]
and
\[
\{u\in {\mathcal B}_2:d(u,{\mathcal T}(h))<\varepsilon \}\subset U_2.
\]
Choose a natural number $n$ with $\frac{\|f-g\|}{2^n}<\varepsilon$. 
Let 
\[
h_k=\frac{k}{2^n}(g-f)+f
\]
for each $0\le k\le 2^n$. By the first part of the proof we have
\[
{\mathcal T}(h_k)+{\mathcal T}(h_{k+2})-2{\mathcal T}(h_{k+1})=0
\qquad \text{($k$)}
\]
holds for $0\le k\le 2^n-2$. For $0\le k\le 2^n-4$,
adding the equations ($k$), 2 times of ($k+1$), and ($k+2$) we have
\[
{\mathcal T}(h_k)+{\mathcal T}(h_{k+4})-2{\mathcal T}(h_{k+2})=0,
\]
whence the equality 
\[
{\mathcal T}(\frac{f+g}{2})=\frac{{\mathcal T}(f)+{\mathcal T}(g)}{2}
\]
holds by induction on $n$.
\end{proof}
Note that an isometry between open sets of Banach algebras need not be 
extended to a linear isometry between these Banach algebras.
\begin{example}\label{hoo}
Let $X=\{x,y\}$ be a compact Hausdorff space consisting of two points.
Let 
\[
{\mathcal U}=\{f\in C(X):\|f\|<1\}\cup
\{f\in C(X):\|f-f_0\|<1\},
\]
where $f_0\in C(X)$ is defined as $f_0(x)=0,\,f_0(y)=10$.
Suppose that 
\[
{\mathcal T}:{\mathcal U}\to {\mathcal U}
\]
is defined as ${\mathcal T}(f)=\tilde f$ if $\|f\|<1$ and 
${\mathcal T}(f)=f$ if $\|f-f_0\|<1$, where 
\begin{equation*}
\tilde f(t)=
\begin{cases}
-f(t),& t=x \\
f(t),& t=y.
\end{cases}
\end{equation*}
Then ${\mathcal T}$ is an isometry from ${\mathcal U}$ onto itself, 
while it cannot be extended to a real linear isometry up to translation.
\end{example}

\section{Linear extensions of isometries}
\label{lin}
In this section we show that isometries between open subgroups of 
the invertibel groups of unital Banach algebras are 
extended to real-linear isometries up to translation.
\begin{lemma}\label{rad}
Let $B$ be a unital Banach algebra and $a\in B$. 
Suppose that $r(ba)=0$ for every $b\in B^{-1}_e$. Then 
$a\in \rad (B)$.
\end{lemma}
\begin{proof}
Suppose that $\alpha$ is a complex number and $0\le t\le 1$. 
Then $r(t\alpha a)=0$, so $-1$ is not in the spectrum of 
$t\alpha a$; $t\alpha a +e_B\in B^{-1}$. 
Thus $\alpha a+e_B$ is in the component of $e_B$; $\alpha
a+e\in B^{-1}_e$.

We will show that $a\in L$ whenever $L$ is a maximal left 
ideal of $B$, which will follow that $a\in \rad (B)$. 
Suppose that there is a maximal left ideal $L$ with 
$a\not\in L$. 
Then $L+Ba$ is a left ideal of $B$ which properly contains $L$, 
so $L+Ba=B$ for $L$ is a maximal left ideal. 
Thus there is $a' \in B$ with $a'a+e_B \in L$. 
Let $\alpha$ be a complex number such that 
$a'-\alpha e_B\in B^{-1}_e$. Since $\alpha a+e_B\in B^{-1}_e$, 
\begin{multline*}
(\alpha a+e_B)^{-1}(a'-\alpha e_B)a+e_B=\\
(\alpha a + e_B)^{-1}(a'-\alpha e_B)a+(\alpha a+e_B)^{-1}
(\alpha a + e_B)\\
=
(\alpha a +e_B)^{-1}(a'a+e_B)\in L
\end{multline*}
hold. 
Thus  $(\alpha a+e_B)^{-1}(a'-\alpha e_B)a+e_B$ is 
singular, hence 
$r((\alpha a+e_B)^{-1}(a'-\alpha e_B)a)>0$, which is  a
contradiction since 
$(\alpha a +e_B)^{-1}(a'- \alpha e_B)\in B^{-1}_e$.
\end{proof}

\begin{theorem}\label{main}
Let $A$ and $B$ be unital Banach algebras and ${\mathfrak A}$ 
and ${\mathfrak B}$ open subgroups of $A$ and $B$ respectively. 
Suppose that $T$ is a surjective isometry from ${\mathfrak A}$ 
onto ${\mathfrak B}$. 
Then there exist a surjective real-linear isometry 
$\tilde T_0$ from $A$ 
onto $B$ and $u_0\in \rad (B)$ which satisfy that 
$T(a)=\tilde T_0(a)+u_0$ for every $a\in {\mathfrak A}$.
\end{theorem}
\begin{proof}
A routine argument shows the existence of 
$\lim_{{\mathfrak A}\ni a \to 0}T(a)$ for an isometry $T$,
 and let $u_0=\lim_{{\mathfrak A}\ni a \to 0}T(a)$. 
Suppose that $f\in B_e^{-1}$ and 
$\sigma (fu_0)\ni \lambda \ne 0$. 
Let $c_{\lambda}=T^{-1}(-\lambda f^{-1})$. It is well-defined 
since $B^{-1}_e\subset {\mathfrak B}$ for ${\mathfrak B}$ 
is an open subgroup of $B^{-1}$. For every $0<s<1$ and 
$0\le t\le 1$, $(1-t)(1-s)c_{\lambda}+tsc_{\lambda} 
\in {\mathfrak A}$ holds, so
\[
T\left(\frac{c_{\lambda}}{2}\right)=
T\left(
\frac{(1-s)c_{\lambda}+sc_{\lambda}}{2}\right)
=\frac{T((1-s)c_{\lambda})+T(sc_{\lambda})}{2}
\]
holds by Lemma \ref{lmu}, and letting $s\to 0$
\[
T\left( 
\frac{c_{\lambda}}{2}
\right)=
\frac
{T(c_{\lambda})+u_0}{2}=
\frac{-\lambda f^{-1}+ u_0}{2},
\]
so 
\[
-\lambda e_B+fu_0=2fT\left(\frac{c_{\lambda}}{2}\right)
\in {\mathfrak B}\subset B^{-1}
\]
holds. Thus $\lambda \not\in \sigma (fu_0)$ holds, 
which is a contradiction proving that 
$\sigma (fu_0)=\{0\}$ or $r(fu_0)=0$ for every 
$f\in B_e^{-1}$. It follows by Lemma \ref{rad} that 
$u_0\in \rad (B)$.

Since $\rad (B)$ is the intersection of the maximal left 
ideal $b^{-1}u\in \rad (B)$ for every $b\in {\mathfrak B}$
and $u\in \rad (B)$. Hence $b^{-1}u+e_B\in B_e^{-1}$, 
so $u+b=b(b^{-1}u+e_B)\in {\mathfrak B}$ holds. 
In particular we see that ${\mathfrak B}-u_0
={\mathfrak B}$. 
Let $T_0:{\mathfrak A}\to \mfB$ be defined as 
$T_0(b)=T(b)-u_0$. Then $T_0$ is well-defined and 
a surjective isometry. Let
$\tilde T_0:A\to B$ be defined as 
\begin{equation*}
\tilde T_0(a)=
\begin{cases}
0, &a=0 \\
T_0(a+2\|a\|e_A)-T_0(2\|a\|e_A),& a\ne 0.
\end{cases}
\end{equation*}
For a non-zero $a\in A$, $a+2\|a\|e_A$ and $2\|a\|e_A$ 
are in $A_e^{-1}$ since 
\[
\|a+2\|a\|e_A-2\|a\|e_A\|=\|a\|<2\|a\|,
\]
and $A^{-1}_e\subset \mfA$ for $\mfA$ is open subgroup of $A$, 
$\tilde T_0$ is well-defined map from $A$ into $B$. 
In the rest of the proof  we show that $\tilde T_0$ is as desired.

We will show that 
\begin{equation}\label{4.0}
T_0(-a)=-T_0(a)
\end{equation}
for every $a\in \mfA$. Let $a\in \mfA$. 
Since $-a\in \mfA$ and $\mfA$ is open, 
$-a+\frac{i}{n}a\in \mfA$ for a sufficiently large 
positive integer $n$. Then by Lemma \ref{lmu} we see 
that 
\[
T_0\left(
\frac{i}{2n}a\right)=
T_0\left(
\frac{a+(-a+\frac{i}{n}a)}{2}\right)
=
\frac{T_0(a)+T_0(-a+\frac{i}{n})}{2}.
\]
Thus 
\[
0=\lim_{n\to \infty}T_0\left(\frac{i}{2n}a\right)
=\lim_{n\to \infty}
\frac{T_0(a)+T_0(-a+\frac{i}{n}a)}{2}=T_0(a)+T_0(-a)
\]
and $T_0(-a)=-T_0(a)$ holds for every $a\in \mfA$. 

By Lemma \ref{lmu}
\[
T_0\left(
\frac{2a+\varepsilon a}{2}\right)
=\frac{T_0(2a)+T_0(\varepsilon a)}{2}
\]
holds for every positive $\varepsilon$ and $a \in \mfA$, 
and letting $\varepsilon
\to 0$, we have $2T_0(a)=T_0(2a)$. 
Applying Lemma \ref{lmu} and by induction on $n$ we see that 
$nT_0(a)=T_0(na)$ holds. 
Substituting $a$ by $\frac{m}{n}a$, 
$T\left(\frac{m}{n}a\right)=
\frac{m}{n}
T_0(a)$
holds for every positive integer $m$ and $n$, and $a\in \mfA$. 
Since $T_0$ is continuous we conclude that 
\begin{equation}\label{homo}
T_0(ra)=rT_0(a)
\end{equation}
holds for every $a\in \mfA$ and $r>0$. 
Then by (\ref{4.0}), (\ref{homo}) holds for every $a\in \mfA$ and 
every real number $r$.

Let $a,b\in \mfA$. Suppose that  
$(1-t)a+tb\in \mfA$ holds for every 
$0\le t\le 1$. Then by Lemma \ref{lmu} 
$T_0\left(\frac{a+b}{2}\right)=
\frac{T_0(a)+T_0(b)}{2}$, 
hence by (\ref{homo})
\begin{equation}\label{sum}
T_0(a+b)=T_0(a)+T_0(b)
\end{equation}
holds. In particular, (\ref{sum}) holds for every pair 
$a$ and $b$ in 
\[
\Omega_A=\{a\in A:\text{$\|a-re_A\|<r$ for some $r>0$}\}
\]
since $\Omega_A$ is convex and 
$\Omega_A\subset A^{-1}_e\subset \mfA$. 
Hence we see that 
\begin{equation}\label{2.5}
\tilde T_0(a)=T_0(a+2\|a\|e_A)-T_0(2\|a\|e_A)=T_0(a)
\end{equation}
holds for every $a\in \Omega_A$ since 
$a+2\|a\|e_A,2\|a\|e_A\in \Omega_A$.

We will show that $\tilde T_0$ is real-linear. Let 
$c$ be a non-zero element in $A$. Then $c+re_A,\,re_A\in 
\Omega_A$ for every positive real number $r\ge 2\|c\|$. 
Hence by (\ref{sum}) 
\begin{multline*}
T_0(c+re_A)+T_0(2\|c\|e_A) \\
=T_0(c+re_A+2\|c\|e_A)=
T_0(c+2\|c\|e_A)+T_0(re_A)
\end{multline*}
holds. Consequently
\begin{equation}\label{3}
\tilde T_0(c)=T_0(c+re_A)-T_0(re_A)
\end{equation}
holds for every $c\in A\setminus \{0\}$ and positive real 
number $r\ge 2\|c\|$. 
We will show that 
$\tilde T_0(a+b)=\tilde T_0(a)+\tilde T_0(b)$ holds 
for every $a,b\in A$. 
If $a$ or $b=0$, the equation is trivial. 
Assume $a\ne 0$ and $b\ne 0$. Since 
$2\|a+b\|\le 2\|a\|+2\|b\|$ we obtain from 
(\ref{3}) that 
\begin{equation}\label{4}
\tilde T_0(a+b)=
T_0(a+b+2\|a\|e_A+2\|b\|e_A)-
T_0(2\|a\|e_A+2\|b\|e_A).
\end{equation}
Applying (\ref{sum}) we observe that 
\[
T_0(a+b+2\|a\|e_A+2\|b\|e_A)=
T_0(a+2\|a\|e_A)+
T_0(b+2\|b\|e_A)
\]
and 
\[
T_0(2\|a\|e_A+2\|b\|e_A)=
T_0(2\|a\|e_A)+T_0(2\|b\|e_A)
\]
holds and we conclude from (\ref{4}) that 
\begin{multline}\label{add}
\tilde T_0(a+b)\\
=T_0(a+2\|a\|e_A)-T_0(2\|a\|e_A)+T_0(b+2\|b\|e_A)
-T_0(2\|b\|e_A)\\
=\tilde T_0(a)+\tilde T_0(b)
\end{multline}
holds, as claimed. 
We will show next that 
$\tilde T_0(ra)=r\tilde T_0(a)$ holds for every $a\in A$ 
and a real number $r$. If $a=0$ or $r=0$, the equation 
is trivial. Let $0\ne a\in A$ and $r\ne0$. Suppose 
$r>0$. Then applying (\ref{homo}) we obtain the 
required equation. Suppose that $r<0$. Since 
$-a+2\|a\|e_A,\,a+2\|a\|e_A,\,2\|a\|e_A \in \Omega_A$ 
equations
\begin{multline*}
T_0(-a+2\|a\|e_A)+T_0(a+2\|a\|e_A) \\
=T_0(2\|a\|e_A+2\|a\|e_A)
=T_0(2\|a\|e_A)+T_0(2\|a\|e_A)
\end{multline*}
hold by (\ref{sum}), and they imply 
\begin{multline}
\tilde T_0(ra)=-r\left(
T_0(-a+2\|a\|e_A)-T_0(2\|a\|e_A)\right) \\
=
-r\left(
-T_0(a+2\|a\|e_A)+T_0(2\|a\|e_A)\right)
=r\tilde T_0(a).
\end{multline}
Henceforth $\tilde T_0$ is real-linear, as claimed.

We will show that $\tilde T_0$ is an isometry. Suppose that 
$a\in A$. Then $\|\tilde T_0(a)\|=\|a\|$ 
if $a=0$. Also 
\begin{multline*}
\|\tilde T_0(a)\|
=\|T_0(a+2\|a\|e_A)-T_0(2\|a\|e_A)\|\\
=\|T(a+2\|a\|e_A)-T(2\|a\|e_A)\|
=\|a+2\|a\|e_A-2\|a\|e_A\|=\|a\|
\end{multline*}
holds if $a\ne 0$. Hence $\tilde T_0$ is an isometry 
since $\tilde T_0$ is linear. 

We will show that $\tilde T_0$ is surjective. 
Let $a \in B$. Let $r$ be a real number with $\|a\|<r$ 
and $\|(T_0(e_A))^{-1}a\|<r$. Then 
\[
\|(T_0(e_A))^{-1}a+re_B-re_B\|<r
\]
implies that $(T_0(e_A))^{-1}a+re_B\in B^{-1}_e$ and by 
(\ref{homo})
\[
a+T_0(re_A)\in \mfB
\]
holds. On the other hand
\begin{equation*}
\|T_0^{-1}(a+T_0(re_A))-re_A\| 
=\|a+T_0(re_A)-T(re_A)\|<r
\end{equation*}
implies that 
$T_0^{-1}(a+T_0(re_A))\in \Omega_A$. Let 
$f=T_0^{-1}(a+T_0(re_A))-re_A$. 
Applying (\ref{sum})
\begin{multline*}
T_0(f+re_A)+T_0(2\|f\|e_A)\\
=T_0(f+re_A+2\|f\|e_A)
=T_0(f+2\|f\|e_A)+T_0(re_A)
\end{multline*}
hold, hence
\begin{multline*}
a=T_0(f+re_A)-T_0(re_A)
=T_0(f+2\|f\|e_A)-T_0(2\|f\|e_A)=
\tilde T_0(f)
\end{multline*}
hold. Consequently $\tilde T_0$ is surjective 
since $a\in B$ is arbitrary.

We will show that $\tilde T_0$ is an extension of 
$T_0$; $\tilde T_0 (a)=T_0(a)$ for every $a\in \mfA$. 
Let $P=\tilde T_0^{-1}\circ T_0:\mfA \to A$. 
Then by (\ref{2.5})
\begin{equation}\label{5}
P(a+2\|a\|e_A)=a+2\|a\|e_A
\end{equation}
holds for every $a\in \mfA$. Since $\tilde T_0^{-1}$ 
is linear
\begin{equation}\label{6}
P(a-2\|a\|e_A)=
-P(-a+2\|a\|e_A)=a-2\|a\|e_A
\end{equation}
holds by (\ref{4.0}) and (\ref{5}) for 
every $a\in \mfA$. For every $0\le t\le 1$
\begin{multline*}
\|t(a+2\|a\|e_A)+(1-t)(\pm 2i\|a\|e_A)
-2(t\pm (1-t)i)\|a\|e_A\| \\
=\|ta\|<2|t\pm (1-t)i|\|a\|
\end{multline*}
holds, for $a\in \mfA$, hence 
\[
t(a+2\|a\|e_A)+(1-t)(\pm2i\|a\|e_A)\in \Omega_A.
\]
Therefore
\begin{multline}
T_0(a+2\|a\|e_A)+T_0(\pm 2i\|a\|e_A)
=
T_0(a+2\|a\|e_A\pm 2i \|a\|e_A)
\end{multline}
holds by (\ref{sum}). Since 
\begin{multline}
\|t(a\pm 2i\|a\|e_A)+(1-t)2\|a\|e_A-
2(\pm it +(1-t))\|a\|e_A\| \\
=\|ta\|<2|\pm it-(1-t)|\|a\|
\end{multline}
holds for every $0\le t \le 1$ and $a\in \mfA$,
\[
T_0(a\pm 2i\|a\|e_A)+T_0(2\|a\|e_A)=
T_0(a\pm 2i\|a\|e_A+2\|a\|e_A)
\]
holds in a same way. It follows that 
\begin{multline}\label{9.3}
\tilde T_0(a)=T_0(a+2\|a\|e_A)-T_0(2\|a\|e_A) \\
=T_0(a\pm 2i\|a\|e_A)-T_0(\pm2i\|a\|e_A)
\end{multline}
holds for every $a\in \mfA$. Applying (\ref{sum})
\[
T_0(\pm 2i\|a\|e_A+4\|a\|e_A)=
T_0(\pm 2i\|a\|e_A)+T_0(4\|a\|e_A)
\]
holds since 
\[
t(\pm 2i\|a\|e_A)+(1-t)4\|a\|e_A\in \mfA
\]
for every $a\in \mfA$ and $0\le t \le 1$. 
Then we obtain
\begin{multline}\label{9.4}
T_0(\pm 2i\|a\|e_A)=T_0(\pm 2i\|a\|e_A+4\|a\|e_A)
-T_0(4\|a\|e_A) \\
=\tilde T_0(\pm 2i\|a\|e_A)
\end{multline}
for every $a\in \mfA$.
It follows by (\ref{9.3}) and (\ref{9.4}) that 
\begin{multline}\label{9.5}
P(a\pm 2i\|a\|e_A)=\tilde T_0^{-1}(\tilde T_0(a)
+T_0(\pm 2i\|a\|e_A)) \\
=\tilde T_0^{-1}(\tilde T_0(a)+\tilde T_0(\pm 2i\|a\|e_A))
=a\pm 2i\|a\|e_A
\end{multline}
holds for every $a\in \mfA$ since $\tilde T_0^{-1}$ is 
linear.

Applying (\ref{5}) and (\ref{6})
\begin{multline}\label{100}
2\|a\|=\|a\pm 2\|a\|e_A-a\|=\|P(a\pm 2\|a\|e_A)-P(a)\| \\
=\|a\pm 2\|a\|e_A -P(a)\|=
\|P(a)-a\pm 2\|a\|e_A\|
\end{multline}
holds for every $a\in \mfA$ 
since $T_0$ and $(\tilde T_0)^{-1}$ are isometric.
In a same way, applying (\ref{9.5})
\begin{equation}\label{100i}
2\|a\|=\|P(a)-a\pm 2i\|a\|e_A\|
\end{equation}
holds for every $a\in \mfA$.

For an element $b\in B$ the numerical range of $b$ is denoted by 
$W(b)$. 
By (\ref{100i}) and \cite[Lemma 2.6.3]{palmer} 
\begin{multline}
\sup \{\mathrm{Im}(\lambda):\lambda \in W(P(a)-a)\} \\
=\inf_{t>0}t^{-1}(\|e_A-it(P(a)-a)\|-1) \\
\le 
2\|a\|(\|e_A-\frac{i}{2\|a\|}(P(a)-a)\|-1)=0.
\end{multline}
Since $W(-P(a)+a)=-W(P(a)-a)$ we have
\begin{multline}
-\inf \{\mathrm{Im} (\lambda):\lambda \in W(P(a)-a)\} \\
=\sup \{\mathrm{Im}(\lambda): \lambda \in W(-P(a)+a)\} \\
=
\inf _{t>0} t^{-1}(\|e_A-it(-P(a)+a)\|-1) \\
\le
2\|a\|(\|e_A+\frac{i}{2\|a\|}(P(a)-a)\|-1)=0
\end{multline}
Thus we see that 
\begin{equation}\label{w1}
W(P(a)-a)\subset \mathbb{R},
\end{equation}
where $\mathbb{R}$ denotes the set of real numbers. 
Applying (\ref{100}) and \cite[Lemma 2.6.3]{palmer} 
in a same way we see that 
\begin{equation}\label{w2}
iW(P(a)-a)=W(i(P(a)-a))\subset \mathbb{R}.
\end{equation}
It follows by (\ref{w1}) and (\ref{w2}) that 
\[
W(P(a)-a)=\{0\}.
\]
Since 
\[
\|P(a)-a\|\le e\|P(a)-a\|_W
\]
holds (
cf. \cite[Theorem 2.6.4]{palmer}), where $\|\cdot \|_W$ denotes the 
numerical radius, 
we see that $P(a)=a$ holds for every $a\in \mfA$.

Since $T(a)=T_0(a)+u_0$ for $a\in \mfA$ by the 
definition of $T_0$, 
we conclude that $T(a)=\tilde T_0(a)+u_0$ holds 
for every $a\in \mfA$.
\end{proof}
\section{Isometrical isomorphisms as metrizable groups}
\label{met}
There are unital commutative $C^*$-algebras $C(X)$ and 
$C(Y)$ which are not real-algebraically isomorphic to each 
other while the invertible groups of these 
algebras are homeomorphically isomorphic to each other
(cf.\cite[Remark 1.7.8]{z2}). 
We show that if the invertible groups are isometrically 
isomorphic to each other, then the corresponding Banach 
algebras are isometrically isomorphic as 
real Banach algebras to each other.
\begin{cor}
Let $A$ and $B$ be unital Banach algebras and 
${\mathfrak A}$ and ${\mathfrak B}$ open subgroups of 
$A^{-1}$ and $B^{-1}$ respectively. Suppose that 
$T$ is a isometrical isomorphism from the metrizable group 
${\mathfrak A}$ onto the metrizable group ${\mathfrak B}$. 
Then $T$ is extended to an isometrical real-algebra 
isomorphism from $A$ onto $B$. In particular, 
$A^{-1}$ is isometrically isomorphic to $B^{-1}$ as a 
metrizable group.
\end{cor}
\begin{proof}
Since $T$ is a group isomorphism $T(e_A)=e_B$ holds. 
Let $a\in {\mathfrak A}$ be an element such that 
$T(a)=2e_B$ and ${\mathfrak A}\ni a_n\to 0$. 
Then 
\[
u_0=\lim_{n\to \infty}T(aa_n)=\lim_{n\to \infty}T(a)T(a_n)
=2e_Bu_0,
\]
where $u_0$ is a radical element described in Theorem 
\ref{main}, hence $u_0=0$. Thus $T$ is extended to a surjective  
real-linear isometry $\tilde T$ from $A$ onto $B$ by Theorem 
\ref{main}. We will show that $\tilde T$ is multiplicative 
Suppose $a,b\in A$. Let $r$ be a sufficiently large real number 
such that $a-re_A,\,b-re_A\in A^{-1}_e$. Since $T$ is 
multiplicative on ${\mathfrak A}$ and $A^{-1}_e\subset 
{\mathfrak A}$,
\[
T((a-re_A)(b-re_A))=T(a-re_A)T(b-re_A)
\]
holds and by a simple calculation in the both side of the 
equation implies that 
\[
\tilde T(ab)=\tilde T(a)\tilde T(b)
\]
holds. Thus $\tilde T$ is a surjective isometrical 
real-algebra isomorphism as desired. Hence $\tilde T(A^{-1})
=B^{-1}$ holds, so $\tilde T$ defines an 
isometrical isomorphism as metrizable groups 
between $A^{-1}$ and $B^{-1}$.
\end{proof}
Applying ${\mathfrak A}=A^{-1}_e$ and 
${\mathfrak B}=B^{-1}_e$ for the above corollary 
the following holds.
\begin{cor}
Let $A$ and $B$ be unital Banach algebras. 
Then $A^{-1}_e$ is isometrically isomorphic to 
$B^{-1}_e$ as a metrizable group if and only if 
$A^{-1}$ is isometrically isomorphic to $B^{-1}$ as a metrizable 
group. In this case $A$ is 
isometrically isomorphic to $B$ as a real Banach algebra.
\end{cor}
\section{Multiplicativity or unti-multiplicativity 
of isometries}
\label{mul}
We proved the following in \cite{hatori}. The proof involves much about 
commutativity and  semisimplicity of the given Banach algebra $A$.
\begin{theorem}\label{st}
Let $A$ be a unital semisimple commutative Banach algebra and $B$ 
a unital Banach algebra. Suppose ${\mathfrak A}$ and ${\mathfrak B}$
are open subgroups of $A^{-1}$ and $B^{-1}$ respectively. 
Suppose that $T$ is a surjective isometry 
from ${\mathfrak A}$ onto
${\mathfrak B}$. Then $B$ is a semisimple and commutative, and 
$(T(e_A))^{-1}T$ is extended to an isometrical real algebra 
isomorphism from $A$ onto $B$. In particular, 
$A^{-1}$ is isometrically isomorphic to 
$B^{-1}$ as a metrizable group.
\end{theorem}
In the following comparison result as the above we make use of 
Theorem \ref{main}.
\begin{cor}\label{comsem}
Let $A$ be a unital commutative Banach algebra and 
$B$ a unital semisimple Banach algebra.
Let ${\mathfrak A}$ be an open subgroup of $A^{-1}$ 
and ${\mathfrak B}$ an open subgroup of $B^{-1}$. 
Suppose that $T$ is a surjective isometry from $\mfA$ onto 
$\mfB$. Then $(T(e_A))^{-1}T$ is extended to a surjective isometrical 
real algebra isomorphism from $A$ onto $B$. Moreover, 
$A$ is semisimple and $B$ is commutative. In particular, 
$A^{-1}$ is isometrically isomorphic to $B^{-1}$ as a metrizable group.
\end{cor}
\begin{proof}
By Theorem \ref{main} there is a $u_0\in \rad (B)$ such that $T-u_0$ is 
extended to a real-linear isometry $\tilde T_0$ 
from $A$ onto $B$. Since 
$B$ is semisimple, $u_0=0$ (cf. \cite[Theorem 2.5.8]{dales}), 
hence  $\tilde T_0$ is an extension of $T$ itself.
We will 
show that $A$ is semisimple. Let $a\in \rad (A)$ and let
$T_a:\mfA\to \mfB$ be defined as $T_a(b)=T(a+b)$ for 
$b\in \mfA$. 
Then $T_a$ is well-defined and a surjective isometry since 
$a+\mfA=\mfA$ 
for $a\in \rad (A)$. By Theorem \ref{main} 
$T_a$ is also extended to a surjective real-linear isometry 
$\tilde T_a$ from 
$A$ onto $B$. For every positive integer
\[
\tilde T_a\left(\frac{e_A}{n}\right)=
T_a\left(\frac{e_A}{n}\right)=T
\left(a+\frac{e_A}{n}\right)=
\tilde T_0\left(a+\frac{e_A}{n}\right)
\]
holds. Letting $n\to \infty$ we have
\[
0=\tilde T_0(a),
\]
hence $a=0$ for $\tilde T_0$ is injective. 
It hollows that $\rad (A)=\{0\}$, 
or $A$ is semisimple. 
Then by Theorem \ref{st} the conclusion holds.
\end{proof}
The hypothesis that $B$ is semisimple in 
Corollary \ref{comsem} is essential 
as the following example (cf. \cite{hatori}) 
shows that a unital isometry 
from $A^{-1}$ onto $B^{-1}$ need not be 
multiplicative nor antimultiplicative 
unless at least one of $A$ or $B$ are  semisimple.
\begin{example}\label{dame}
Let 
\[
A_0=\{
\left(
\begin{smallmatrix}
0&a&b \\0&0&c \\ 0&0&0
\end{smallmatrix}
\right)
:a,\,\,b,\,\,c\in {\mathbb C}\}.
\]
Let 
\[
A=\{
\left(
\begin{smallmatrix}
\alpha&a&b \\0&\alpha&c \\ 0&0&\alpha
\end{smallmatrix}
\right)
:\alpha,\,\,a,\,\,b,\,\,c\in {\mathbb C}\}
\]
be the unitization of $A_0$, where the multiplication (in $A_0$) is the 
zero multiplication; $MN=0$ for every $M,N\in A_0$.
Let $B=A$ as sets, while the multiplication in $B$ is the usual multiplication 
for matrices. 
Then $A$ and $B$ are unital Banach algebras under the usual operator norm. 
Note that $A$ is commutative and $A$ nor $B$ are not semisimple. 
Note also that 
$A^{-1}=\{\left(
\begin{smallmatrix}
\alpha&a&b \\0&\alpha&c \\ 0&0&\alpha
\end{smallmatrix}
\right)\in A:\alpha\ne 0\}$ and $B^{-1}=\{\left(
\begin{smallmatrix}
\alpha&a&b \\0&\alpha&c \\ 0&0&\alpha
\end{smallmatrix}
\right)\in B:\alpha\ne 0\}$.
Define 
$T:A^{-1}\to B^{-1}$ by $T(M)=M$. 
Then $T$ is well-defined and a surjective unital 
($T(e_A)=e_B$) isometry.
On the other hand $A^{-1}$ is not (group) isomorphic to $B^{-1}$, 
in particular, $T$ is not multiplicative nor antimultiplicative.
\end{example}
We show a result for standard operator algebras.
\begin{cor}\label{standard}
Let ${\mathcal X}$ (resp. ${\mathcal Y}$) 
be a Banach space. Suppose that $A$ (resp. $B$) 
is a unital closed 
subalgebra of ${\mathcal B}({\mathcal X})$ 
(resp. ${\mathcal B}({\mathcal Y})$), 
the Banach algebra of all the bounded operators on 
${\mathcal X}$ (resp. ${\mathcal Y}$), 
which contains all finite rank operators. 
Suppose that $T$ is a surjective isometry from 
an open subgroup ${\mathfrak A}$ of $A^{-1}$ 
onto an open subgroup $\mfB$ of 
$B^{-1}$. Then there exists an invertible bounded linear 
or conjugate linear operator 
$U:X\to Y$ such that $T(a)=T(e_A)UaU^{-1}$ 
for every 
$a\in \mfA$, or there exists an invertible 
bounded linear or conjugate 
linear operator $V:X^*\to Y$ such that 
$T(a)=T(e_A)Va^*V^{-1}$ for every $a\in \mfA$. 
In particular, if $T$ is unital in the sense that 
$T(e_A)=e_B$, then $T$ is multiplicative 
(in this case, $A^{-1}$is isometrically
isomorphic to $B^{-1}$ as a metrizable group) 
or anti-multiplicative (in this case, $A^{-1}$is isometrically
anti-isomorphic to $B^{-1}$ as a metrizable group).
\end{cor}
\begin{proof}
By Theorem \ref{main} there is a $u_0\in \rad (B)$ 
such that $T-u_0$ is 
extended to a real-linear isometry 
$\tilde T_0$ from $A$ onto $B$. Since 
$B$ is semisimple (cf. \cite[Theorem 2.5.8]{dales}), 
$u_0=0$ follows, hence 
$\tilde T_0=T$ on $\mfA$. 
Thus $(T(e_A))^{-1}\tilde T_0$ is an
additive surjection such that 
$(T(e_A))^{-1}\tilde T_0(A^{-1})
=B^{-1}$. Applying  Theorem 3.2 in \cite{houcui} for 
$(T(e_A))^{-1}\tilde T_0$,
there exists an invertible bounded linear or 
conjugate linear operator 
$U:{\mathcal X}\to {\mathcal Y}$ such that 
$(T(e_A))^{-1}\tilde T_0(a)=UaU^{-1}$ ($a\in A$), 
or there exists an invertible bounded linear or conjugate linear 
operator 
$V:{\mathcal X}^*\to {\mathcal Y}$ 
such that $(T(e_A))^{-1}\tilde T_0(a)=Va^*V^{-1}$ ($a\in A$). 
Henceforth the conclusions hold.
\end{proof}
Let $M_n$ be the algebra of all $n\times n$ matrices over the complex number 
field. For $M\in M_n$ the spectrum  is denoted by $\sigma (M)$ and $M^t$ 
is the transpose of $M$. $E$ denotes the identity matrix. 
Let $\|\cdot\|$ and $\|\cdot\|'$ 
denote any matrix norms on $M_n$ (cf. \cite{hojo}).
\begin{cor}
If $S$ is a surjection from an open subgroup 
$\mfA$ of the group $M_n^{-1}$ of the invertible $n\times n$ 
matrices over the complex number field onto 
an open subgroup $\mfB$ of $M_n^{-1}$ such that 
$\|S({M})-S({N})\|'=\|{M}-{N}\|$ for all ${M},{N}\in \mfA$, 
then 
there exists an invertible matrix ${U}\in M_n$ such that 
$S({M})=S(E){U}{M}{U}^{-1}$ for all ${M}\in \mfA$, or
$S({M})=S(E){U}{M}^t{U}^{-1}$ for all ${M}\in \mfA$, or 
$S(M)=S(E)U\overline{M}U^{-1}$  for all $M\in \mfA$, or
$S(M)=S(E)U\overline{M}^tU^{-1}$  for all $M\in \mfA$ hold.
In particular, if $S$ is unital, then $S$ is multiplicative or 
anti-multiplicative.
\end{cor}
\begin{proof}
By Corollary \ref{standard} there is an invertible matrix $U$ such that 
one of 
the following four occurs.
\begin{enumerate}
\item
$S(M)=S(E)UMU^{-1}$ holds for every $M\in M_n^{-1}$, 
\item
$S(M)=S(E)U\overline{M}U^{-1}$ holds for every $M\in M_n^{-1}$, 
\item
$S(M)=S(E)UM^tU^{-1}$ holds for every $M\in M_n^{-1}$, 
\item
$S(M)=S(E)U\overline{M}^tU^{-1}$ holds for every $M\in M_n^{-1}$. 
\end{enumerate}
Henceforth the conclusion holds.
\end{proof}



\begin{thebibliography}{99}



\bibitem{b}
S.~Banach,
"Theory of Linear Operations", North-Holland, Amsterdam, New York, 
Oxford, Tokyo 1987



\bibitem{dales}
H.~G.~Dales,
"Banach Algebras and Automatic Continuity", 
London Mathematical Society Monographs. New Series, 24. Oxford Science Publications. The Clarendon Press, Oxford University Press, New York, 2000






\bibitem{hatori}
O.~Hatori,
{\it Isometries between groups of invertible elements in Banach algebras},
to appear


\bibitem{hojo}
R.~A.~Horn and C.~R.~Johnson,
"Matrix analysis
 Corrected reprint of the 1985 original", Cambridge University Press, Cambridge, 1990

\bibitem{houcui}
J.~Hou and J.~Cui,
{\it Additive maps on standard operator algebras preserving invertibilities or zero divisors},
Linear Algebra Appl.,
{\bf 359}(2003), 219--233



\bibitem{ja1}
K.~Jarosz,
{\it The uniqueness of multiplication in function algebras},
Proc. Amer. Math. Soc.,
{\bf 89}(1983), 249--253































\bibitem{mu}
S.~Mazur and S.~Ulam,
{\it Sur les transformationes isom\'etriques 
d'espaces vectoriels norm\'es}, 
C. R. Acad. Sci, Paris,
{\bf 194}(1932), 946--948






\bibitem{palmer}
T.~W.~Palmer,
"Banach algebras and the general theory of $\sp *$-algebras. Vol. I. Algebras and Banach algebras",
Encyclopedia of Mathematics and its Applications, 49. Cambridge University Press, Cambridge, 1994.






\bibitem{v}
J.~V\"ais\"al\"a, 
{\it A proof of the Mazur-Ulam theorem}, 
Amer. Math. Monthly, 
{\bf 110}(7)(2003), 633--635



\bibitem{z}
W.~\. Zelazko,
{\it A characterization of multiplicative linear functionals 
in complex Banach algebras}, 
Studia Math.,
{\bf 30}(1968), 83--85

\bibitem{z2}
W.~\. Zelazko, 
"Banach Algebras", Elsivier, 1973


















\end{thebibliography}
\end{document}